\documentclass[11pt,twoside]{article}%{atmp}

\usepackage{amsmath,amssymb}

\usepackage{color}

\usepackage{latexsym,amssymb,epsfig,amsmath,amsthm}
\usepackage{graphicx}
\usepackage[all]{xy}
\usepackage[english]{babel}
\usepackage{url}

\newcommand{\IZ}{\mathbb{Z}}
\newcommand{\IC}{\mathbb{C}}
\newcommand{\IF}{\mathbb{F}}
\newcommand{\IP}{\mathbb{P}}
\newcommand\cH{{\mathcal H}}
\newcommand{\ba}{\begin{array}}
\newcommand{\ea}{\end{array}}

\newtheorem{theorem}{Theorem}[section]
\newtheorem{corollary}[theorem]{Corollary}
\newtheorem{lemma}[theorem]{Lemma}
\newtheorem{proposition}[theorem]{Proposition}
\theoremstyle{definition}
\newtheorem{definition}[theorem]{Definition}
\theoremstyle{remark}

\newtheorem{example}[theorem]{Example}

\begin{document}

\sloppy

\title{Toroidal Orbifolds \`a la Vafa-Witten}
\author{Jimmy Dillies}
%\address{Department of Mathematics - Texas A\&M University\\Mailstop 3368\\College Station, TX 77843-3368} 

\maketitle

%\addressemail{jimmy@math.tamu.edu}

\begin{abstract}
We classify orbifolds obtained by taking the quotient of a three torus by Abelian extensions of  $\IZ /n\times \IZ /n$ automorphisms, where each torus has a multiplicative $\IZ/n$ action ($n\in\{3,4,6\}$).
This 'completes' the classification of orbifolds of the above type initiated by Donagi and Faraggi (\cite{DoFa}) and, Donagi and Wendland (\cite{DoWe}).
\end{abstract}

\section{Introduction}
%%%%%%%%%%%%%%%%%%%%%%

In 1985, Dixon, Harvey, Vafa and Witten pioneered the study of string theory on orbifolds (\cite{DHVW2}). A recurrent model, the $\IZ /2\times \IZ /2$ orbifold, which was introduced by Vafa and Witten (\cite{VaWi}) has been studied extensively. In particular, Donagi and Faraggi classified further quotients using symmetric shifts (\cite{DoFa}). They deduce that the three generation vacua are not obtained in this manner. Seeking a better model, Donagi and Wendland (\cite{DoWe}) studied quotients of three tori by Abelian extensions of the $\IZ /2\times \IZ /2$ automorphisms. \\
We classify here the orbifolds obtained by taking the quotient of three tori by Abelian extensions of  $\IZ /n\times \IZ /n$ automorphisms, where each torus has a multiplicative $\IZ/n$ action ($n\in\{3,4,6\}$). 

\subsection*{Acknowledgments}
%%%%%%%%%%%%%%%%%%%%%%%%%%%%%%
This work is part of the author's thesis at the University of Pennsylvania which has been a benevolent host. It was partially supported through a fellowship of the BAEF (Belgian American Educational Foundation).\\
I would like to thank Professor Ron Donagi for suggesting this problem and for many fruitful conversations.

%\cutpage %move this line so that the first page breaks at the appropriate place.

%\setcounter{page}{<insert page # for second page>}

\noindent

%%%%%%%%%%%%%%%%%%%%%
\section{Construction and Results}
%%%%%%%%%%%%%%%%%%%%%%

%%%%%%%%%%%%%%%%%%
\subsection{Construction}
%%%%%%%%%%%%%%%%%%

Consider the elliptic curve $E_3$, quotient of the complex plane by the sublattice $\Lambda_n=\IZ \oplus \omega_3 \IZ$ where $\omega_3=e^\frac{i\pi}{3}$. This curve is endowed with a multiplicative automorphism: if we write $[x]$ for the class in $E_3$ of $x\in\IC$ and $\zeta_3$ for a primitive third root of unity, the map is given explicitly by $[x]\mapsto[\zeta_3 x]$. This automorphism generates an action of $\IZ/3\IZ$ on $E_3$.
Consequently, if we take the variety obtained by taking three copies of $E_3$, $X_3=E_3\times E_3\times E_3$, it comes with an action of $(\IZ/3\IZ)^3$. This action restricts to an action of $(\IZ/3\IZ)^2$ as schematized in the following diagram

\centerline{
\xymatrix{
0 \ar[r] & (\IZ/3\IZ)^2 \ar[r] \ar@{^{(}-->}[dr] &  (\IZ/3\IZ)^3 \ar[r]^{+} \ar@{^{(}->}[d] & \IZ/3\IZ \ar[r] & 0 \\
 & & Aut(X_3) & & 
}
}

Define $T_3$ to be the set of points of $X_3$ which are fixed by $(\IZ/3\IZ)^2$. It is a subgroup of $E_3[3]\times E_3[3]\times E_3[3]$ isomorphic to $(\IZ/3\IZ)^3$ ($E_3[3]$ denotes the 3-torsion points of the elliptic curve). The group $T_3$ acts by translation on $X_3$ and its action commutes with the one of $(\IZ/3\IZ)^2$ which was just introduced. We have therefore an action of the direct product $V_3:=(\IZ/3\IZ)^2\times(\IZ/3\IZ)^3$ which we will call the Vafa-Witten action of order $3$. \\
In the following sections we will deal with quotients of $X_3$ by subgroups of $V_3$ which project onto the multiplicative $(\IZ/3\IZ)^2$.\\
The above construction can be mimiced with the elliptic curves $E_4=\IC/(\IZ\oplus i\IZ)$ and $E_6=E_3$ respectively endowed with the automorphisms $[x]\mapsto[\zeta_6 x]$ and $[x]\mapsto[\zeta_4 x]$. As a result, we construct an action of $V_4=((\IZ/4\IZ)^2\times (\IZ/2\IZ)^3)$ and $V_6=((\IZ/6\IZ)^2\times \textrm{Id})$ respectively on $X_4$ and $X_6$; these actions will be called the Vafa-Witten actions of order $4$ and $6$. 
As for $X_3$ we will be interested in quotients of $X_{n=4,6}$ by subgroups of $V_{n=4,6}$ which project onto the multiplicative part.

\noindent We can synthesize the situation in the following table:
$$
\ba{cccc}
n  & \omega_n & \zeta_n & T_n \\
\hline
\hline
3 & e^\frac{i\pi}{3} & e^\frac{2i\pi}{3} & (\IZ/3\IZ)^3 \subset E_3[3]^3\\
4 & i & i & (\IZ/2\IZ)^3 \subset E_4[2]^3 \\
6 & e^\frac{i\pi}{3} & e^\frac{i\pi}{3} & \{0\} \\
%\hline
\ea
$$

%%%%%%%%%%%%%%%%%%%%%
\subsection{Notation}
%%%%%%%%%%%%%%%%%%%%%

The element $g:([z_1],[z_2],[z_3])\mapsto ([\zeta_i^{m_1}z_1+a_1 t_n],[\zeta_i^{m_2}z_2+a_2 t_n],[\zeta_i^{m_3}z_3+a_3 t_n])$
will be denoted $g=(m_1,m_2,m_3;a_1,a_2,a_3)$. 
The element $t_n$ is a generator of $T_n|_{E_n}$. Notice that the $m_i$ add up to zero in $\IZ/n\IZ$.\\
We will call the $m_i$'s the {\it twist} part, and the $a_i$'s the {\it shift} part. 
\begin{definition}
The \textit{rank} of $G$ is its number of generators minus $2$. 
\end{definition}
The definition is such that, when $G$ is a direct product of the multiplicative group by some subgroup of $T_n$, the rank corresponds to the number of generators of the translation part.

%%%%%%%%%%%%%%%%%%%%%
\subsection{Results}
%%%%%%%%%%%%%%%%%%%%%

\begin{definition}
We call $\cH_n$ the set of orbifolds obtained by taking the quotient of $X_n$ by subgroups of the $n^\textrm{th}$ Vafa-Witten group, $V_n$, which surject onto the multiplicative component. 
\end{definition}

In the next sections we get the following classification:

\begin{proposition}
\label{resultch1}
The sets $\mathcal H_n$ contain only finitely many homeomorphism classes of orbifolds. The number of classes is $8$ for $\mathcal H_3$ and $\mathcal H_4$, and is $1$ for $\mathcal H_6$. 
\end{proposition}

The Hodge numbers of the spaces are given in the tables which follow. For each homeomorphism class of $\cH_n$, we give a representative group $G$ by listing its generators, as well as the Hodge numbers of $X_n/G$.

\subsubsection*{Case $n=3$}
%%%%%%%%%%%%%%%%%%%%%%%%%%%

For $n=3$, the generators are (where applicable): $(1,2,0;a_1,a_2,a_3)$, $(2,0,1;b_1,b_2,b_3)$, $(0,0,0;c_1,c_2,c_3)$ and $(0,0,0;d_1,d_2,d_3)$.
$$
\begin{array}{ccccccc}
%\hline 
 \# & (a_1,a_2,a_3) & (b_1,b_2,b_3)  & (c_1,c_2,c_3) & (d_1,d_2,d_3) & (h_{11},h_{12}) \\
\hline
\hline 
 III.1 & (0,0,0) & (0,0,0) &&& (84,0) \\
 III.2 & (0,0,0) & (0,1,0) &&& (24,12) \\
 III.3 & (0,0,0) & (1,1,0) &&& (18,6) \\
 III.4 & (0,0,1) & (1,1,0) &&& (12,0)  \\
\hline 
 III.5 & (0,0,0) & (0,0,0) & (0,1,1) && (40,4) \\ 
 III.6 & "& "& (1,1,1) && (36,0) \\
% III.7 & (0,0,0) & (0,1,0) & (1,1,0) && (40,4)\\
 III.7 & (0,0,0) & (0,1,0) & (1,0,1) && (16,4) \\
 III.8 & "& "& (1,1,1) && (18,6) \\
\end{array}
$$

All the above orbifolds are simply connected except $III.4$ whose fundamental group is $\IZ/3\IZ$.

\subsubsection*{Case $n=4$}
%%%%%%%%%%%%%%%%%%%%%%%%%%%

 For $n=4$, the generators are (where applicable): $(1,3,0;a_1,a_2,a_3)$, $(3,0,1;b_1,b_2,b_3)$, $(0,0,0;c_1,c_2,c_3)$ and $(0,0,0;d_1,d_2,d_3)$.

$$
\begin{array}{cccccc}
\# & (a,b,c) & (a',b',c') & (c_1,c_2,c_3) & (d_1,d_2,d_3) & (h_{11},h_{12}) \\
\hline
\hline
IV.1 & (0,0,0) & (0,0,0) &&& (90,0) \\
IV.2 & (0,0,0) & (0,1,0) &&& (54,0) \\
IV.3 & (0,0,0) & (1,1,0) &&& (42,0) \\
IV.4 & (0,0,1) & (1,1,0) &&& (30,0) \\
\hline
 IV.5 & (0,0,0) & (0,0,0) & (0,1,1) && (61,1) \\
 IV.6 & "& "& (1,1,1) && (54,0) \\
% IV.7 & (0,0,0) & (0,1,0) & (1,1,0) && (61,1) \\
 IV.7 & (0,0,0)  & (0,1,0) & (1,0,1) && (38,0)\\
 IV.8 & "& "& (1,1,1) && (42,0)\\
\end{array}
$$

All these orbifolds have trivial fundamental group.

\subsubsection*{Case $n=6$}
%%%%%%%%%%%%%%%%%%%%%%%%%%%

For $n=6$, there is a unique case corresponding to the quotient by the  Vafa-Witten group. The orbifold $X_6/V_6$ has Hodge numbers $(80,0)$ and is simply connected.

%%%%%%%%%%%%%%%%%%%%%%%%%%%%%%%%%%%%%%%%%%%%%%%
\section{Homeomorphism classes of $\mathcal H_n$}
%%%%%%%%%%%%%%%%%%%%%%%%%%%%%%%%%%%%%%%%%%%%%%%

In this section we will classify the elements of $\mathcal H_n$ up to homeomorphism. The classification will be made according to the rank of the group acting on $X_n$. 

%%%%%%%%%%%%%%%%%%%%%%%%%%%
\subsection{General Lemmas}
%%%%%%%%%%%%%%%%%%%%%%%%%%%%

To identify groups, we will make recurrent use of the following lemma:

\begin{lemma}
Let $G$ be a subgroup of $V_4$ of $V_4$ (resp. $V_3$) which surjects onto the multiplicative part, then $G$ has at least two generators. The first two generators can be taken of the form $g_1=(1,2,0;*,*,*)$ and $g_2=(2,0,1;*,*,*)$ (resp. $g_1=(1,3,0;*,*,*)$ and $g_2=(3,0,1;*,*,*)$). If there are more than two generators, then they can be taken of the form $g_{i>2}=(0,0,0;*,*,*)$.
\end{lemma}

\begin{proof}
Since $G$ must surject onto the multiplicative part, it clearly has at least two generators.  These two generators $g_1,g_2$ can be chosen to be the lifts of the generators of the multiplicative group. For $n=4$ (resp. $n=3$) this group admits as generators $(1,2,0)$ and $(2,0,1)$ (resp.  $(1,3,0)$ and $(3,0,1)$), so the nature of the first two generators is settled.\\
Let $g_i$ be another generator of $G$. Since the twist part of $g_1$ and $g_2$ generates the multiplicative part, there exists a word $w$ in the first two generators so that $w.g_i=(0,0,0;*,*,*)$. We can now substitute $g_i$ by $w.g_i$. 
\end{proof}

The most powerful tool to identify groups will be conjugation by torsion elements of $E_n ^3$.
For simplicity we will work on a single torus. Consider the transformation $z \mapsto \zeta_i^a z+\tau t_i$. We will conjugate it with the translation by $\lambda$, any given element of $E_n$:
$$(\zeta_i^a(z+\lambda)+\tau t_i)-\lambda=\zeta_i^a z + \underbrace{\tau t_i + (\zeta_i^a \lambda -\lambda)}_\textrm{new shift}$$  
If $a\neq 0$, we can choose a $\lambda$ such that $\tau t_i + (\zeta_i^a \lambda -\lambda)=0$, which implies that our transformation is conjugate to $\zeta_i^a z$.

%%%%%%%%%%%%%%%%%%%%%%%%%%%%%%%
\subsection{Quotients of $X_3$}
%%%%%%%%%%%%%%%%%%%%%%%%%%%%%%%

%%%%%%%%%%%%%%%%%%%%%%%
\subsubsection*{rank 0}
%%%%%%%%%%%%%%%%%%%%%%%

Although there are a priori $(3^3)^2=729$ possible choices of $g_1=(1,2,0;a_1,a_2,a_3)$ and $g_2=(2,0,1;b_1,b_2,b_3)$, by using conjugation and symmetry, we can, to begin with, restrict ourselves to at most $4$ cases:

\begin{lemma}
Without loss of generality we can assume $g_1$ and $g_2$ to be of the form $(1,2,0;0,0,a_3)$ and $(2,0,1;b_1,b_2,0)$ with $a_3,b_1,b_2\in \{0,1\}$.
\end{lemma}

\begin{proof}
Let $G$ be a group with generators $g_1=(1,2,0;a_1,a_2,a_3)$ and 
$g_2=(2,0,1;b_1,b_2,b_3)$. Conjugate both elements by $(0,0,0;\alpha_1,\alpha_2,\alpha_3)$. The element $g_1$ gets mapped to 
$(1,2,0;a_1+(\zeta \alpha_1 -\alpha_1),a_2+(\zeta^2 \alpha_2 -\alpha_2),a_3)$ while $g_2$ gets mapped to $(2,0,1;b_1+(\zeta^2 \alpha_1 -\alpha_1),b_2,b_3+(\zeta \alpha_3 -\alpha_3))$. We can take $\alpha_{1\ldots 3}$ so that  
$$a_1+(\zeta \alpha_1 -\alpha_1)=a_2+(\zeta^2 \alpha_2 -\alpha_2)=b_3+(\zeta \alpha_3 -\alpha_3))=0$$
Finally, the symmetry between $t_3$ and $2t_3$ allows us to take the remaining entries in $\{0,1\}$.
\end{proof}

\begin{lemma}
If we take generators as in the previous lemma, the entries $b_1$ and $b_2$ are symmetric.
\end{lemma}

\begin{proof}
Note that $g_1=(1,2,0;0,0,\delta_3)$ with $\delta_3\in \{0,1\}$ and $g_2=(2,0,1;b_2,b_3,0)$.
The group spanned by $g_1,g_2$ is the same as the group spanned by $g_1^2,g_1 g_2$, that is
$(2,1,0;0,0,2\delta_3)$ and $(0,2,1;b_1,b_2,\delta_3)$. We now rearrange the order of the tori: $(1\;2\;3)\leadsto(2\;1\;3)$ so that the new generators read $(1,2,0;0,0,\delta_3)$ and $(2,0,1;b_2,b_1,\delta_3)$. By conjugating by an appropriate element of the third torus we get as second generator the required $(2,0,1;b_2,b_1,0)$. 
\end{proof}

The above lemmas restrict the number of cases to 6:

\begin{itemize}
\item $g_1=(1,2,0;0,0,1)$ $g_2=(2,0,1;0,0,0)$
\item $g_1=(1,2,0;0,0,1)$ $g_2=(2,0,1;1,0,0)$
\item $g_1=(1,2,0;0,0,1)$ $g_2=(2,0,1;1,1,0)$
\item $g_1=(1,2,0;0,0,0)$ $g_2=(2,0,1;0,0,0)$
\item $g_1=(1,2,0;0,0,0)$ $g_2=(2,0,1;1,0,0)$
\item $g_1=(1,2,0;0,0,0)$ $g_2=(2,0,1;1,1,0)$
\end{itemize}

We will show that among those $6$ classes, $2$ are redundant.
Also, we will simplify the notation further: e.g. $\underline{a}=(1,1,1)$ will denote the element $g_1=(1,2,0;1,1,1)$ while $\underline{b}=(1,1,0)$ will denote the element $g_2=(2,0,1;1,1,0)$.

\begin{lemma}
The group generated by $\underline{a}=(0,0,1)$ and $\underline{b}=(0,0,0)$ is isomorphic to the group spanned by $\underline{a}=(0,0,0)$, $\underline{b}=(0,1,0)$.
\end{lemma}

\begin{proof}
We can replace the generators $g_1$ and $g_2$ by their squares: $(2,1,0;0,0,2)$ and $(1,0,2;0,0,0)$. If we rearrange the terms in the order $(1\;2\;3)\leadsto(1\;3\;2)$ and we permute the generators we get $(1,2,0;0,0,0)$ and $(2,0,1;0,2,0)$. This is what we want up to relabeling. 
\end{proof}

\begin{lemma}
The elements $\underline{a}=(0,0,1)$ and $\underline{b}=(1,0,0)$ and the elements $\underline{a}=(0,0,0)$, $\underline{b}=(1,1,0)$ generate isomorphic groups:
\end{lemma}

\begin{proof}
The basis $g_2,g_1^2 g_2^2$ is equivalent to $g_1,g_2$. It is made out of the vectors 
$(2,0,1;1,0,0)$ and $(0,1,2;2,0,2)$. We now rearrange the tori using the permutation $(1\;2\;3)\leadsto(3\;1\;2)$ to get the basis $(1,2,0;0,1,0)$,$(2,0,1;2,2,0)$. We now conjugate with an appropriate translation on the second tori to get $(1,2,0;0,0,0)$,$(2,0,1;2,2,0)$. 
\end{proof}

As a conclusion we have, 
\begin{proposition}
There are four homeomorphism classes of quotients of $X_3$ by groups of rank $0$ in $\mathcal H_3$. We have written a representative of each class in the following table:
$$
\begin{array}{ccccc}
%\hline 
 \# & (a_1,a_2,a_3) & (b_1,b_2,b_3)  & (h_{11},h_{12}) & \pi_1 \\
\hline
\hline 
 III.1 & (0,0,0) & (0,0,0) & (84,0) & 1 \\
 III.2 & (0,0,0) & (0,1,0) & (24,12) & 1 \\
 III.3 & (0,0,0) & (1,1,0) & (18,6) & 1\\
 III.4 & (0,0,1) & (1,1,0) & (12,0) & \IZ/3 \\
%\hline
\end{array}
$$
\end{proposition}

\begin{proof}
We have seen through the previous lemmas that there are at most $4$ types of isomorphism of groups. By computing the Hodge diamond of the associated Calabi-Yau threefold (see section 3) we deduce that they yield  four different varieties.
\end{proof}

%%%%%%%%%%%%%%%%%%%%%%
\subsubsection*{rank 1}
%%%%%%%%%%%%%%%%%%%%%%

We label the third generator $g_3=(0,0,0;c_1,c_2,c_3)$ or $\underline{c}$.
We will extend the list of rank $0$ groups using the following rules:

\begin{lemma} (Donagi \& Faraggi \cite{DoFa}) \textit{(Reduction principle)}
Let $G$ be a group of which one of the generators is of the form $(0,0,0;x_1,x_2,x_3)$ and exactly one of the $x_i\neq0$. The quotient $X_n/G$ is then homeomorphic to $X_n/\bar{G}$ where $\bar{G}$ is the quotient of $G$ by the subgroup spanned by that generator.
\end{lemma}

The idea is that if there is an element which consists in a translation on a unique curve, we can first take the quotient by this element and have another three torus on which the rest of the group acts. \\
As a corollary, we can restrict ourselves to $g_3$'s where at least two of the $c_k$'s are not zero. Furthermore, we have the following simplifications:

\begin{lemma}
If $c_k\neq 0$ and $a_k=b_k=0$ we can assume $c_k=1$.
\end{lemma}

\begin{proof}
It follows from the symmetry between $t_3$ and $2t_3$.
\end{proof}

\begin{lemma}
We can assume that the shift part of $\underline{c}$ is not a non-zero multiple of the shift part of $\underline{a}$ or of $\underline{b}$.
\end{lemma}

\begin{proof}
Without loss of generality, assume that $\underline{c}$ is a multiple of $\underline{a}$, then, we can substitute $g_1$ by $g_1 g_3^k$ (so the shift part is $(0,0,0)$) to get a new first generator without translation. In other words, we have reduced the group to a previous case.  
\end{proof}

We will now try to discern the groups:

\begin{enumerate}

\item {$\underline{a}=(0,0,0);\underline{b}=(0,0,0)$}\\
We can choose $\underline{c}$ to be either $(1,1,0)$ or $(1,1,1)$. All other cases resume to these two using the previous points and $S_3$ symmetry.

\item {$\underline{a}=(0,0,0);\underline{b}=(0,1,0)$}\\ 
We can choose, using the previous points, a $\underline{c}$ of the form $(\delta_1,c_2,\delta_3)$ where the $\delta_i\in\{0,1\}$. \\
However, we can also assume that $c_2\in\{0,1\}$:

\begin{proof}
The generator $g_3$ is equivalent to $g_3^2$, so since we can relabel the first and third translations without loss of generality, we see that we can assume $c_2$ to be in $\{0,1\}$.  
\end{proof}

It now seems that we have $4$ possible choices for $\underline{c}$, namely $(1,1,0)$, $(0,1,1)$, $(1,0,1)$ and $(1,1,1)$. We will show that the first two are actually redundant:\\

\begin{proof}
Adjoin $(0,0,0;0,1,0)$ to the group generated by . $\underline{a}=(0,0,0);\underline{b}=(0,1,0)$ and $\underline{c}=(1,1,0)$. It is easy to see that up to $S_3$ symmetry this group is equivalent to the action of the group obtained by joining $(0,0,0;0,1,0)$ to the group generated by $\underline{a}=(0,0,0);\underline{b}=(0,0,0)$ and $\underline{c}=(0,1,1)$. By the reduction principle, all these groups generate the same space up to homeomorphism, so the first case reduces to an ancient case.\\
For the second case, notice that the group generated by $g_1,g_2,g_3$ is the same as the group spanned by $g_1,g_2 g_3^2, g_3$. The element $g_2 g_3^2=(2,0,1;0,0,2)$. If we conjugate these new elements by an adequate translation on the third torus, the first and last generator are unchanged while $g_2 g_3^2\leadsto (2,0,1;0,0,0)$.  
\end{proof}

The computation of the Hodge numbers will assure us that the $2$ remaining cases are independent.

\item {$\underline{a}=(0,0,0);\underline{b}=(1,1,0)$}\\
We can choose $\underline{c}$ to be of the form $(c_1,c_2,\delta_3)$ where $\delta_3\in\{0,1\}$ (same argument as previously). Now, we could  replace $g_3$ by its square and, since we can change the last component of $\underline{c}$ freely, it means that we could replace $(c_1,c_2,\delta_3)$ by $(c_1^2,c_2^2,\delta_3)$. We will do this so to have a minimal number of entries equal to $2$ in $(c_1,c_2)$.\\
Using the above rules, the possible $\underline{c}$'s are $(1,2,0)$,$(0,1,1)$,$(1,0,1)$,$(1,1,1)$ or $(1,2,1)$.\\
Actually, none of these cases is new:

\begin{itemize}
\item $(1,2,0)$: We have $(g_1,g_2,g_3)=(g_1,g_2 g_3^2,g_3)$. The element $g_2 g_3^2=(2,0,1;0,2,0)$. After conjugation, and up to transforming the third generator $g_3$, we can let $g_2 g_3^2\leadsto(2,0,1;0,0,0)$.

\item $(0,1,1)$: We have $(g_1,g_2,g_3)=(g_1^2,g_1 g_2 g_3^2,g_3)$, where $g_1^2=(2,1,0;0,0,0)$ and $g_1 g_2 g_3^2=(0,2,1;1,0,2)$. Conjugating by some appropriate translation on the third torus, the second generator becomes\\
 $g_1 g_2 g_3^2\leadsto(0,2,1;1,0,0)$. We now reorder the tori $(1,2,3)\leadsto(2,1,3)$ and we get a previous case.

\item $(1,0,1)$: We have $(g_1,g_2,g_3)=(g_1,g_2 g_3^2,g_3)$. The element $g_2 g_3^2=(2,0,1;0,1,2)$.
After conjugating with an element of translation on the third torus, the second generator can be taken
to be $(2,0,1;0,1,0)$.

\item $(1,1,1)$: We have $(g_1,g_2,g_3)=(g_1,g_2 g_3^2,g_3)$. The element $g_2 g_3^2=(2,0,1;0,0,1)$.
Again, we can conjugate by an element of translation on the third torus, so to have our second generator $(2,0,1;0,0,0)$. 

\item $(1,2,1)$:  We have $(g_1,g_2,g_3)=(g_1,g_2 g_3^2,g_3)$. The element $g_2 g_3^2=(2,0,1;0,2,2)$.
We can conjugate by an element of translation on the third torus, so to have our second generator $\leadsto(2,0,1;0,2,0)$. Up to renaming, we again reduced to a previous case. 

\end{itemize}

So we conclude that there are no interesting extensions in this case.

\item {$\underline{a}=(0,0,1);\underline{b}=(1,1,0)$}\\ 
We can choose $\underline{c}$ to be of the form $(c_1,c_2,c_3)$. Let us first undercover some symmetry:
we have $(g_1,g_2,g_3)=(g_1^2,g_1 g_2,g_3)$ and if we permute the order of the tori $(1,2,3)\leadsto(2,1,3)$ we get the generators $(1,2,0;0,0,2)$,$(2,0,1;1,1,1)$ and $(0,0,0;c_2,c_1,c_3)$. Now we can conjugate by a translation on the third torus to let the second generator become $(2,0,1;1,1,0)$ and leave the other two unchanged. Now we can relabel the new $a_3$ into a $1$, and we must therefore substitute the new $c_3$ by its square. Finally, we get the generators
$(1,2,0;0,0,1)$, $(2,0,1;1,1,0)$ and $(0,0,0;c_2,c_1,c_3^2)$.\\
So we can let $c_3$ be $0$ or $1$ up to a permutation of the $c_1,c_2$. \\
Using the above symmetry, we have the following possibilities which we show to be reducible to previous cases: 

\begin{itemize}
\item $(1,2,0)$: We have $(g_1,g_2,g_3)=(g_1 ,g_2 g_3, g_3)$ where $g_2 g_3=(2,0,1;2,0,0)$. Using the symmetry of the first translation entries of the second generator (up to variation of $g_3$) we have reduced to a previous case.
\item $(1,0,1)$: We have $(g_1,g_2,g_3)=(g_1 g_3^2,g_2,g_3)$ where $g_1 g_3^2=(1,2,0;2,0,0)$. Up to changing $g_2$ we can conjugate by a translation element on the first torus to get the first generator to become $(1,2,0;0,0,0)$. So we reduced to a previous case.
\item $(0,1,1)$: We have $(g_1,g_2,g_3)=(g_1 g_3^2,g_2,g_3)$ where $g_1 g_3^2=(1,2,0;0,2,0)$. We can conjugate by a translation element on the second torus to get the first generator to become $(1,2,0;0,0,0)$. So we reduced to a previous case.
\item $(1,1,1)$: We have $(g_1,g_2,g_3)=(g_1,g_2 g_3^2,g_3)$ where $g_2 g_3^2=(2,0,1;0,0,2)$. We can conjugate by a translation element on the third torus to get the second generator to become $(2,0,1;0,0,0)$.
\item $(2,1,1)$: We have $(g_1,g_2,g_3)=(g_1 g_3^2,g_2,g_3)$ where $g_1 g_3^2=(1,2,0;1,2,0)$. Up to changing $g_2$ we can conjugate by a translation element on the first and second torus to get the first generator to become $(1,2,0;0,0,0)$. So we reduced to a previous case.
\item $(1,2,1)$: We have $(g_1,g_2,g_3)=(g_1 g_3^2,g_2,g_3)$ where $g_1 g_3^2=(1,2,0;1,2,0)$. Up to changing $g_2$ we can conjugate by a translation element on the first and second torus to get the first generator to become $(1,2,0;0,0,0)$. So we reduced to a previous case.
\end{itemize}

\end{enumerate}

From the above discussion we conclude:

\begin{proposition}
There are four homeomorphism classes in $\cH_3$ coming from groups of rank $1$. We have written a representative of each class in the following table:
$$
\begin{array}{cccccc}
% \hline
 \# & (a_1,a_2,a_3) & (b_1,b_2,b_3) & (c_1,c_2,c_3) & (h_{11},h_{12}) & \pi_1 \\
\hline 
\hline 
 III.5 & (0,0,0) & (0,0,0) & (0,1,1) & (40,4) & 1\\ 
 III.6 & & & (1,1,1) & (36,0) &  1\\
% III.7 & (0,0,0) & (0,1,0) & (1,1,0) & (40,4)& 1\\
 III.7 & (0,0,0) & (0,1,0) & (1,0,1) & (16,4) & 1\\
 III.8 & & & (1,1,1) & (18,6) & 1\\
\end{array}
$$
\end{proposition}

%%%%%%%%%%%%%%%%%%%%%%
\subsubsection*{rank 2}
%%%%%%%%%%%%%%%%%%%%%%

Since we do not consider rank two groups which we can reduce to a lower rank, we have restrictions on the third and fourth generator: by the reduction principle they cannot generate an element with only one non-zero entry in the shift part.\\
Consider the third and fourth generators as elements of $\IF_3^3$ (they just have shift parts). 

\begin{lemma}
There are exactly 4 linear planes in $\IF_3^3$ which do not intersect the coordinate axes outside the origin.
\end{lemma}

\begin{figure}
\centering
\includegraphics[height=3cm]{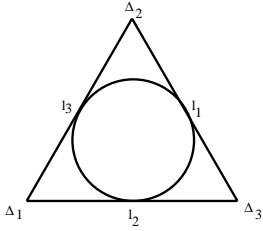}\caption{Planes in $\IF_3^3$ not crossing the coordinate axes outside the origin}
\label{f33}
\end{figure}

\begin{proof}
Let $H$ be a plane which does not intersect the coordinate axes outside the origin. Since it contains the origin, it intersects the three coordinate planes in a line. For each coordinate plane $P_i$, the only two possible lines are the diagonal $\Delta_i$ and, the only other line which is not a coordinate axis, $l_i$. 
The choice of any two lines. not in the same coordinate plane, out of $\{l_1,l_2,\Delta_1,\Delta_2\}$ gives an adequate plane. In particular, there are 4 of them. If we include $\Delta_3$ and $l_3$ in the picture, we look at the coplanarity of the lines, which is readily checked and can be visualized in figure \ref{f33}. Each plane is represented by either an edge of the triangle or the inscribed circle.
\end{proof}

Remarks:
\begin{itemize}
\item The three planes represented by the edges of the triangle in figure \ref{f33} are clearly $S_3$ symmetric. Therefore, it is enough to consider the adjunction of any one of those three to the first case of rank 0 subgroups.

\item Since all rank 1 cases associated to $\underline{a}=(0,0,0),\underline{b}=(1,1,0)$ and  $\underline{a}=(0,0,1),\underline{b}=(1,1,0)$ were reduced to previous cases, all rank 2 extensions will also be reducible, so we just need to deal with $\underline{a}=(0,0,0),\underline{b}=(0,0,0)$
and $\underline{a}=(0,0,0),\underline{b}=(0,1,0)$.
\end{itemize}

We will now work out the rank 2 cases:
\begin{itemize}
\item Using $S_3$ symmetry, for $\underline{a}=\underline{b}=(0,0,0)$, we have as possible extension $\underline{c}=(0,1,1),\underline{d}=(1,1,0)$ and $\underline{c}=(0,1,2),\underline{d}=(1,2,0)$. Now, we can replace $\underline{d}$ by its square and since $c_3$ and $d_1$ are freely chosen, we see that this case is equivalent to the previous one.\\
 
\item All $4$ possible extensions of $\underline{a}=(0,0,0),\underline{b}=(0,1,0)$ contain either the element $(0,2,2)$ or $(0,2,1)$. Using this element it is easy to see that we can reduce $g_2$ to $(2,0,1;0,0,0)$: multiply $g_2$ by this element and conjugate by an appropriate translation on the third coordinate. We show hereby that there are no further rank 2 cases.
\end{itemize}

A priori, there is thus only one rank 2 case:

$$
\begin{array}{cccc}
 (a_1,a_2,a_3) & (b_1,b_2,b_3) & (c_1,c_2,c_3) & (d_1,d_2,d_3) \\
\hline 
\hline
 (0,0,0) & (0,0,0) & (0,1,1) & (1,1,0) \\ 
\end{array}
$$
However, it can be reduced to the $III.1$ case : by applying the reduction principle of Donagi and Faraggi, this case is equivalent to the Vafa-Witten group case (adjoin any translation) which in turn is equivalent to the $III.1$ case.

%%%%%%%%%%%%%%%%%%%%%%
\subsubsection*{rank 3}
%%%%%%%%%%%%%%%%%%%%%%

The only group is the Vafa-Witten group. It reduces to the III.1 case.

%%%%%%%%%%%%%%%%%%%%%%%%%%%%%%
\subsection{Quotients of $X_4$}
%%%%%%%%%%%%%%%%%%%%%%%%%%%%%%%

Since there are only two fixed points per torus, the structure of the translation locus is simpler. A translation element will be of the form $(\delta_1,\delta_2,\delta_3)$ with $\delta_i\in\{0,1\}$. 

%%%%%%%%%%%%%%%%%%%%%%
\subsubsection*{rank 0}
%%%%%%%%%%%%%%%%%%%%%%

Using the same argument as for $\IZ/3\IZ$ (being careful to replace square by inverse) we quickly get the following list:

\begin{proposition}
There are four homeomorphism classes in $\mathcal H_3$ obtained by taking the quotients by rank $0$ groups. We have written a representative of each class in the following table:
$$
\begin{array}{cccc}
\# & (a,b,c) & (a',b',c')  & (h_{11},h_{12}) \\
\hline
\hline
IV.1 & (0,0,0) & (0,0,0) & (90,0) \\
IV.2 & (0,0,0) & (0,1,0) & (54,0) \\
IV.3 & (0,0,0) & (1,1,0) & (42,0) \\
IV.4 & (0,0,1) & (1,1,0) & (30,0) \\
\end{array}
$$
\end{proposition}

%%%%%%%%%%%%%%%%%%%%%%
\subsubsection*{rank 1}
%%%%%%%%%%%%%%%%%%%%%%

Again, we will use the same arguments as for $n=3$. 

\begin{lemma}
If $\underline{a}=(0,0,0)$, $\underline{b}=(0,0,0)$ we can assume $\underline{c}$ to be of the form
\begin{itemize}
\item $(0,1,1)$.
\item $(1,1,1)$.
\end{itemize} 
\end{lemma}

\begin{proof}
Using $S_3$ symmetry, all other $\underline{c}$'s with two non-zero entries are equivalent to the one listed here.
Also, the reduction principle excludes all $\underline{c}$'s with a single non-zero entry.
\end{proof}

\begin{lemma} 
If $\underline{a}=(0,0,0)$, $\underline{b}=(0,1,0)$ then we can assume that $\underline{c}$ is one of the following: 
\begin{itemize}
%\item $(1,1,0)$. 
\item $(1,0,1)$. 
\item $(1,1,1)$. 
\end{itemize}
\end{lemma}

\begin{proof}
The case where $\underline{c}$ is $(0,1,1)$ is reducible: replace $g_2$ by $g_2 g_3=(3,0,1;0,0,1)$ and conjugate by a translation element on the third torus to get a second generator which is translation free.\\
The case where $\underline{c}$ is $(1,1,0)$ is reducible: adjoin the element $(0,0,0;0,1,0)$ and use the same argument as for the order $3$ case.
\end{proof}

\begin{lemma}
If $\underline{a}=(0,0,0)$, $\underline{b}=(1,1,0)$ then there are no new cases.
\end{lemma}

\begin{proof}
If $\underline{c}$ is
\begin{itemize}
\item $(0,1,1)$. This case is reducible:  we replace $g_1$ by $g_1^{-1}$ and $g_2$ by $g_2 g_3=(0,3,1;1,0,1)$. Secondly we conjugate by a translation element on the third torus to let the second generator become $(0,3,1;1,0,0)$. Finally, we arrange the torus in the order $(1,2,3)\leadsto(2,1,3)$. We reduced to an extension of the previous type of rank 0.
\item $(1,0,1)$. This case is reducible:  replace $g_2$ by $g_2 g_3=(3,0,1;0,1,1)$ and conjugate by a translation element on the third torus to get a second generator of the form $(3,0,1;0,1,0)$. We reduced to an extension of the previous type of rank 0.
\item $(1,1,1)$. This case is reducible: replace $g_2$ by $g_2 g_3=(3,0,1;0,0,1)$ and conjugate by a translation element on the third torus to get a second generator which is translation free.
\end{itemize}
\end{proof}

\begin{lemma}
If $\underline{a}=(0,0,1)$, $\underline{b}=(1,1,0)$ then there are no new cases.
\end{lemma}

\begin{proof}
If $\underline{c}$ is
\begin{itemize}
\item $(0,1,1)$. We have $(g_1,g_2,g_3)=(g_1 g_3,g_2,g_3)$ where $g_1 g_3=(1,3,0;0,1,0)$. We can conjugate by a translation element on the second torus to get the first generator to become $(1,3,0;0,0,0)$. So we reduced to a previous case.
\item $(1,0,1)$. We have $(g_1,g_2,g_3)=(g_1 g_3,g_2,g_3)$ where $g_1 g_3=(1,3,0;1,0,0)$. Up to changing $g_2$ we can conjugate by a translation element on the first torus to get the first generator to become $(1,3,0;0,0,0)$. So we reduced to a previous case. 
\item $(1,1,1)$. We have $(g_1,g_2,g_3)=(g_1 g_3,g_2,g_3)$ where $g_1 g_3=(1,3,0;1,1,0)$. Up to changing $g_2$ we can conjugate by a translation element on the first and second torus to get the first generator to become $(1,3,0;0,0,0)$. So we reduced to a previous case.
\end{itemize}
\end{proof}

Using those lemmas, we conclude 

\begin{proposition}
There are four homeomorphism classes generated by groups of rank $1$ in $\mathcal H_3$. We have written a representative of each class in the following table: 
$$
\begin{array}{ccccc}
\#  & (a_1,a_2,a_3) & (b_1,b_2,b_3) & (c_1,c_2,c_3) & (h_{11},h_{12}) \\
\hline 
\hline
 IV.5 & (0,0,0) & (0,0,0) & (0,1,1) & (61,1) \\
 IV.6 & & & (1,1,1) & (54,0) \\
% IV.7 & (0,0,0) & (0,1,0) & (1,1,0) & (61,1) \\
 IV.7 & (0,0,0)& (0,1,0)& (1,0,1) & (38,0)\\
 IV.8 & & & (1,1,1) & (42,0)\\
%\hline
\end{array}
$$
\end{proposition}

%%%%%%%%%%%%%%%%%%%%%%
\subsubsection*{rank 2}
%%%%%%%%%%%%%%%%%%%%%%

Since we do not want cases which reduce to lower ranks, we need the 3rd and 4th generator to span 
 a subgroup which does not contain elements where only one of the entries is non-zero.  There is actually a unique possibility. In geometric language:
 
\begin{lemma}
There is a unique $2$ dimensional vector subspaces of $\IF_2^3$ which does not intersect the coordinate axes outside the origin.
\end{lemma}

\begin{proof}
Let $H$ be a plane verifying the above conditions. Since $H$ passes through the origin, it must intersect each coordinate plane in at least a line. Since $H$ does not intersect the coordinate axes, the intersection lines must be the first diagonals. Therefore $\{(0,0,0),(0,1,1),(1,0,1),(1,1,0)\}\subset H$ and since a plane over $\IF_2$ contains $4$ elements, we are done. 
\end{proof}

Since the group of translations is isomorphic to $\IF_3^3$, we get the direct result that

\begin{corollary}
There is at most one rank 2 group associated to each group of rank 0 -- which we classified earlier.
\end{corollary}

\begin{corollary}
There is no non-reducible rank 2 case associated to the last three types of rank 0. 
\end{corollary}

\begin{proof}
In the last two types of rank 0, one generator has a shift with 2 non-zero entries. Since this shift belongs to the plane $H$, the case can be reduced to a previous one.\\
In the second case of rank 0, $\underline{b}=(0,1,0)$, so replacing $g_2$ with its composition with $(0,0,0;0,1,1)$ we get  $(3,0,1;0,0,1)$. We can now conjugate with an adequate translation on the third torus to get the second generator in the form $(3,0,1;0,0,0)$, i.e. we are back to the first type of rank 0. 
\end{proof}

So we conclude that apparently there is only one type of rank 2 subgroup, namely:
$$
\begin{array}{cccc}
% \hline
 (a_1,a_2,a_3) & (b_1,b_2,b_3) & (c_1,c_2,c_3) & (d_1,d_2,d_3) \\
\hline
\hline
(0,0,0) & (0,0,0) & (0,1,1) & (1,1,0) \\
%\hline
\end{array}
$$
However, as in the case of $X_3$, it is reducible to $IV.1$.

%%%%%%%%%%%%%%%%%%%%%%
\subsubsection*{rank 3}
%%%%%%%%%%%%%%%%%%%%%%

The only case is the whole group $V_4$. It is reducible to the $III.1$ case.

%%%%%%%%%%%%%%%%%%%%%%%%%%%%%%%
\subsection{Quotients of $V_6$}
%%%%%%%%%%%%%%%%%%%%%%%%%%%%%%%

There is a unique case as there are no translations commuting with the multiplicative $\IZ/6\IZ\times \IZ/6\IZ$ action. The Hodge numbers of the quotient are $(80,0)$.

%%%%%%%%%%%%%%%%%%%%%%%%%%%%%
\section{The Hodge Structure}
%%%%%%%%%%%%%%%%%%%%%%%%%%%%%

To compute the orbifold cohomology of the spaces that we have classified in the previous section, we will use the orbifold cohomology introduced in physics (see e.g. \cite{DHVW2}) and formalized by Ruan \& al. (\cite{CR}):

\begin{definition}
The orbifold cohomology of the quotient of the manifold $X_n$ by the finite group $G$, which is also the cohomology of a crepant resolution of this quotient, is 
$$H^{i,j}(X_n/G)=\bigoplus_{
\begin{array}{c}
[g]\in G\\
U \in X_n ^g
\end{array}
}H^{i-\kappa(g),j-\kappa(g)}(U)^{C(g)}$$
where we sum over the conjugacy classes of $G$ and the components of the fixed locus of each element. $C(g)$ denotes the centralizer of $g$ in $G$.
The set $X_n^g$ is the fixed locus of any element in the conjugacy class of $g$, i.e. $\{x\in X_n  \;\textrm{such that}\; g.x=x\}$; $U$ denotes one of its components. Moreover, if the action of $g$ sends $[z_i]_{i=1\ldots 3}$ to $[e^{2\pi i\theta_i}z_i]_{i=1\ldots 3}$ with $0\leq\theta_i < 1$ then the we define the shift function
 \footnote{Also known as \textit{fermionic shift}, \textit{degree shifting number} or \textit{age}}
 by $\kappa(g)=\sum_{i=1}^3 \theta_i$. \\
\end{definition}
In our case, the definition of the Vafa-Witten groups forces the function $\kappa$ to take its values in $\{0,1,2\}$. 
Also, since we are dealing with Abelian groups and given that the local action is essentially unique, the formula simplifies to
$$H^{i,j}(X_n/G)=\bigoplus_{
\begin{array}{c}
g\in G
\end{array}
}H^{i-\kappa(g),j-\kappa(g)}(X_n^g)^G .$$

Note that the action of an element of $G$ on the cohomology depends uniquely on the multiplicative part: the action of $\zeta_i$ is $\zeta_i.dz_i=\zeta_i dz_i$ and $\zeta_i.d\bar z_i=\bar{\zeta}_i d\bar z_i$ while the action of $T_i$ is trivial. 

\begin{example}
Take $g=(m_1,m_2,m_3;a_1,a_2,a_3)$ and $\Omega=dz_1\wedge dz_2 \wedge d\bar{z}_1 \wedge d\bar{z}_3$ we have 
$$\Omega \stackrel{g}{\longmapsto} \zeta_i^{(m_1+m_2)-(m_1+m_3)}\Omega = \zeta_i^{(m_2-m_3)}\Omega$$
\end{example}

In order to compute the Hodge structure of the $X_n/G$ we first need to classify the possible fixed loci of an element of $g\in G$ on $X_n$. 

\begin{lemma}
Consider the element $g=(m_1,m_2,m_3;a_1,a_2,a_3)\in G$, its fixed locus $X_n ^g$ can be of four different topological types
\end{lemma}
 
\begin{proof}

There are four exclusive forms of $g$ which correspond to the four different fixed loci:

\begin{enumerate}

\item The identity element $(0,0,0;0,0,0)$ has as fixed locus the whole variety $X_n$. 

\item If for a certain index $k\in\{1,2,3\}$ we have that $m_k=0$ and $a_k\neq 0$ then the fixed locus of $g$ is the empty set. Indeed, translation acts fixed point freely on an elliptic curve and therefore the action of $g$ on the whole of $X_n$ is also fixed point free.

\item The fixed locus of $g$ where all $m_k$'s are different from $0$ is a collection of points.

\item The fixed locus of $g$, where for exactly one of the k's $m_k=a_k=0$ (and excluding case 2), is made out of elliptic curves (the number depends on $n$).  

\end{enumerate}

Note that if two of the $m_k$'s are zero then the third one also has to be zero and therefore all cases have been exhausted in the above list. \\

\end{proof}

We will now compute the $G$ invariant cohomology from these fixed loci:

\begin{lemma} For all $n\in\{3,4,6\}$,  the G-invariant part of the cohomology of $X_n$ is\footnote{We have tilted the Hodge diamonds $45^o$ to the left to facilitate typesetting.}

$$H^*(X_n)^G=
\ba{cccc}
1 & 0 & 0 & 1\\
0 & 3 & 0 & 0\\
0 & 0 & 3 & 0\\
1 & 0 & 0 & 1\\
\ea
$$

\end{lemma}

\begin{proof}
The Hodge diamond of the 3-torus is 

$$
\ba{cccc}
1 & 3 & 3 & 1\\
3 & 9 & 9 & 3\\
3 & 9 & 9 & 3\\
1 & 3 & 3 & 1\\
\ea
$$
The invariance of $H^{0,0}$ and $H^{3,0}$ is straightforward. \\
For the other components, note that the action of $G$ on $H^{1,0}$ and $H^{0,1}$, and thus on the whole of $H(X_n)$, is diagonal with respect to the standard basis. Therefore, it will be enough to check the behavior of the basis elements to find the G-invariant part.
Let $dz_k$ be a generator of $H^{1,0}$, it is not fixed by the element $(1,1,n-2;*,*,*)$. 
Therefore we have ${h^{1,0}}^G=0$.\\
Similarly, $dz_k \wedge dz_l$, a generator of  $H^{2,0}$, is not fixed by the element $(1,1,n-2;*,*,*)$ and thus ${h^{2,0}}^G=0$. \\
Only the generators of $H^{1,1}$ of the form $dz_i \wedge d\bar z_i$ are invariant, the others are killed by the element $(1,n-1,0;*,*,*)$. Therefore the dimension of the G-invariant part of have $H^{1,1}$ is $3$.\\
Consider now $H^{2,1}$; by symmetry we can restrict ourselves to the generators $dz_1 \wedge dz_2 \wedge d\bar z_3$ and $dz_1 \wedge dz_2 \wedge\ d\bar z_1$. The former is not fixed by $(1,0,n-1;*,*,*)$, while the latter is not fixed by $(1,1,n-2;*,*,*)$. We conclude that ${h^{2,1}}^G=0$. Finally, the diamond is completed using $dz_k \leftrightarrow d\bar{z}_k$ and Hodge symmetry.

\end{proof}

\begin{lemma} 
Assume that after identification via $G$-action, $g\in G$ has as fixed locus a collection of $n$ points. The contribution
to the cohomology of $g$ and $g^{-1}$, $H^*(X^g)^G\oplus H^*(X^{g^{-1}})^G$ is equivalent to the contribution of $n$ projective lines: $H^*(n\IP^1)$.
\end{lemma}

\begin{proof}
The elements $g$ and $g^{-1}$ have the same fixed locus, whose cohomology is exclusively $H^{0,0}$. Given that the fixed locus of $g$ is made of points, we know that $\kappa(g)$ is non-zero and that none of the $\theta_i$ is 0. We claim that $\{\kappa(g),\kappa(g^{-1})\}$ is exactly $\{1,2\}$. Indeed, if we denote by $\theta'_i$ the linearized action of $g^{-1}$ on the $i^{th}$ component, then we have the relation $\theta'_i=1-\theta_i$. Therefore $\kappa(g^{-1})=3-\kappa(g)$.

\end{proof}

\begin{lemma}
The $G$-invariant part of the cohomology of a fixed elliptic curve will be either the cohomology of the projective line {\tiny $\ba{cc} 1&0\\ 0&1\ea$} or the one of an elliptic curve {\tiny $\ba{cc}1&1\\1&1 \ea$}. In either case $\kappa$ will be 1.
\end{lemma}

\begin{proof}
Since $H^{0,0}$ is obviously invariant, only $H^{1,0}$ which is of dimension 1 might not preserved, in which case the invariant cohomology corresponds to that of a $\IP^1$.\\
Since we do not deal with the fixed locus of the identity, $\kappa$ is not $0$. Moreover, we know that the action is trivial on one component. Therefore we have that $\kappa$, which is the sum $\theta_1+\theta_2+\theta_3$ where one $\theta_i=0$ and the two others of norm less than one, can only be one.

\end{proof}

%%%%%%%%%%%%%%%%%%%%%
\subsection{Notation}
%%%%%%%%%%%%%%%%%%%%%

We are now ready to compute the Hodge numbers of each orbifold.\\
Each group will be represented, as previously, by the shift part of its generators.
We will list the group elements $g$ which have a nonempty fixed locus $X^g$ and their contribution, i.e. for each group we will have a collection $\{(g,(h^{1,1},h^{1,2}))\}$. To lighten the notation we will 

\begin{enumerate}
\item Count twice the contribution of non-trivial elements which are not  involutions but have a union of curves as fixed locus. In compensation we will not be writing their inverse. 
\item If a non involutive element fixes a union of points, then we will count its contribution together with the one of its inverse. We will not write down its inverse.
\end{enumerate}

%%%%%%%%%%%%%%%%%%%%
\subsection{Order 3}
%%%%%%%%%%%%%%%%%%%%

We have eight cases to compute:

\begin{itemize}

\item \fbox{$(0,0,0)(0,0,0) \leadsto (h_{11},h_{12})=(84,0)$}\\
$\{$
($(0,0,0;0,0,0),(3,0)$), 
($(1,2,0;0,0,0),(18,0)$), 
($(2,0,1;0,0,0),(18,0)$),
($(0,1,2;0,0,0),(18,0)$), 
($(1,1,1;0,0,0),(27,0)$) $\}$
\item \fbox{$(0,0,0)(0,1,0)\leadsto (24,12)$}\\
$\{$
($(0,0,0;0,0,0),(3,0)$), 
($(1,2,0;0,0,0),(6,6)$), 
($(0,2,1;0,1,0),(6,6)$), 
($(1,1,1;0,1,0),(9,0)$)$\}$
\item \fbox{$(0,0,0)(1,1,0)\leadsto(18,6)$}\\
$\{$($(0,0,0;0,0,0),(3,0)$), 
($(1,2,0;0,0,0),(6,6)$), 
($(1,1,1;1,1,0),(9,0)$)$\}$
\item \fbox{$(0,0,1)(1,1,0)\leadsto(12,0)$}\\
$\{$
($(0,0,0;0,0,0),(3,0)$), 
($(1,1,1;1,1,2),(9,0)$)$\}$

\item \fbox{$(0,0,0)(0,0,0)(0,1,1)\leadsto(40,4)$}\\
$\{$($(0,0,0;0,0,0),(3,0)$), 
($(1,2,0;0,0,0),(6,0)$), 
($(2,0,1;0,0,0),(6,0)$), 
($(0,1,2;0,0,0),(6,0)$), 
($(0,1,2;0,1,1),(2,2)$), 
($(0,1,2;0,2,2),(2,2)$), 
($(1,1,1;0,0,0),(9,0)$), 
($(1,1,1;0,1,1),(3,0)$),
($(1,1,1;0,2,2),(3,0)$)$\}$

\item \fbox{$(0,0,0)(0,0,0)(1,1,1)\leadsto(36,0)$}\\$\{$(
$(0,0,0;0,0,0),(3,0)$),(
$(1,2,0;0,0,0),(6,0)$),(
$(2,0,1;0,0,0),(6,0)$),(
$(0,1,2;0,0,0),(6,0)$),(
$(1,1,1;0,0,0),(9,0)$),(
$(1,1,1;1,1,1),(3,0)$),(
$(1,1,1;2,2,2),(3,0)$)$\}$

%\item \fbox{$(0,0,0)(0,1,0)(1,1,0)\leadsto(40,4)$}\\$\{$(
%$(0,0,0;0,0,0),(3,0)$),(
%$(1,2,0;0,0,0),(2,2)$),(
%$(1,2,0;1,1,0),(6,0)$),(
%$(1,2,0;2,2,0),(6,0)$),(
%$(2,0,1;2,0,0),(2,2)$),(
%$(0,2,1;0,1,0),(6,0)$),(
%$(1,1,1;0,1,0),(3,0)$),(
%$(1,1,1;1,2,0),(9,0)$),(
%$(1,1,1;2,0,0),(3,0)$)$\}$

\item \fbox{$(0,0,0,)(0,1,0)(1,0,1)\leadsto(16,4)$}\\
$\{$
($(0,0,0;0,0,0),(3,0)$),(
$(1,2,0;0,0,0),(2,2)$),(
$(0,2,1;0,1,0),(2,2)$),(
$(1,1,1;0,1,0),(3,0)$),(
$(1,1,1;1,1,1),(3,0)$),(
$(1,1,1;2,1,2),(3,0)$)
$\}$

\item \fbox{$(0,0,0)(0,1,0)(1,1,1)\leadsto(18,6)$}\\
$\{$
($(0,0,0;0,0,0),(3,0)$),(
$(1,2,0;0,0,0),(2,2)$),(
$(2,0,1;2,0,2),(2,2)$),(
$(0,2,1;0,1,0),(2,2)$),(
$(1,1,1;0,1,0),(3,0)$),(
$(1,1,1;1,2,1),(3,0)$),(
$(1,1,1;2,0,2),(3,0)$)
$\}$

\end{itemize}

%%%%%%%%%%%%%%%%
\subsection{Order 4}
%%%%%%%%%%%%%%%%

For this family, we also have eight cases to compute: 
\begin{itemize}

\item \fbox{$(0,0,0)(0,0,0)\leadsto (h_{11},h_{12})=(90,0)$}\\
$\{$
($(0,0,0;0,0,0),(3,0)$),(
$(1,3,0;0,0,0),(8,0)$),(
$(3,0,1;0,0,0),(8,0)$),(
$(0,3,1;0,0,0),(8,0)$),(
$(1,1,2;0,0,0),(12,0)$),(
$(1,2,1;0,0,0),(12,0)$),(
$(2,1,1;0,0,0),(12,0)$),(
$(2,2,0;0,0,0),(9,0)$),(
$(2,0,2;0,0,0),(9,0)$),(
$(0,2,2;0,0,0),(9,0)$)
$\}$

\item \fbox{$(0,0,0)(1,1,0)\leadsto(42,0)$}\\
$\{$
($(0,0,0;0,0,0),(3,0)$),
($(1,3,0;0,0,0),(4,0)$),
($(1,1,2;0,0,0),(8,0)$), 
($(1,2,1;1,1,0),(4,0)$),
($(2,1,1;1,1,0),(4,0)$),
($(2,2,0;0,0,0),(5,0)$),
($(2,0,2;0,0,0),(7,0)$),
($(0,2,2;0,0,0),(7,0)$)$\}$

\item \fbox{$(0,0,0)(0,1,0)\leadsto(54,0)$}\\
$\{$
($(0,0,0;0,0,0),(3,0)$),
($(1,3,0;0,0,0),(4,0)$),
($(0,1,3;0,1,0),(4,0)$),
($(1,1,2;0,0,0),(8,0)$),
($(1,2,1;0,1,0),(4,0)$),
($(2,1,1;0,1,,[3],0),(8,0)$),  
($(2,2,0;0,0,0),(7,0)$),
($(2,0,2;0,0,0),(9,0)$),
($(0,2,2;0,0,0),(7,0)$) 
$\}$

\item \fbox{$(0,0,1)(1,1,0)\leadsto(30,0)$}\\
$\{$
($(0,0,0;0,0,0),(3,0)$),
($(1,1,2;0,0,0),(4,0)$), 
($(1,2,1;0,1,0),(4,0)$), 
($(2,1,1;0,1,0),(4,0)$), 
($(2,2,0;0,0,0),(5,0)$),
($(2,0,2;0,0,0),(5,0)$),
($(0,2,2;0,0,0),(5,0)$)
$\}$

\item \fbox{$(0,0,0)(0,0,0)(0,1,1)\leadsto(61,1)$}\\
$\{$,
($(0,0,0;0,0,0),(3,0)$)
($(1,3,0;0,0,0),(4,0)$),
($(3,0,1;0,0,0),(4,0)$),
($(0,3,1;0,0,0),(4,0)$),
($(1,1,2;0,0,0),(6,0)$),
($(1,2,1;0,0,0),(6,0)$),
($(2,1,1;0,0,0),(6,0)$),
($(2,2,0;0,0,0),(6,0)$),
($(2,0,2;0,0,0),(6,0)$), 
($(0,2,2;0,0,0),(5,0)$),
($(0,3,1;0,1,1),(2,0)$),
($(1,1,2;0,1,1),(2,0)$),
($(1,2,1;0,1,1),(2,0)$),
($(2,1,1;0,1,1),(4,0)$),
($(0,2,2;0,1,1),(1,1)$)
$\}$

\item \fbox{$(0,0,0)(0,0,0)(1,1,1)\leadsto(54,0)$}\\
$\{$
($(0,0,0;0,0,0)$,$(3,0)$),
($(1,3,0;0,0,0),(4,0)$),
($(3,0,1;0,0,0),(4,0)$),
($(0,3,1;0,0,0),(4,0)$),
($(1,1,2;0,0,0),(6,0)$),
($(1,2,1;0,0,0),(6,0)$),
($(2,1,1;0,0,0),(6,0)$),
($(2,2,0;0,0,0),(5,0)$),
($(2,0,2;0,0,0),(5,0)$),
($(0,2,2;0,0,0),(5,0)$),
($(1,1,2;1,1,1),(2,0)$),
($(1,2,1;1,1,1),(2,0)$),
($(2,1,1;1,1,1),(2,0)$)
$\}$

%\item \fbox{$(0,0,0)(0,1,0)(1,1,0)\leadsto(42,0)$}\\
%$\{$
%($(0,0,0;0,0,0),(3,0)$),
%$(1,3,0;0,0,0),(2,0)$),(
%$(3,0,1;1,0,1),(2,0)$),( 
%$(0,1,3;0,1,0),(2,0)$),(
%$(1,1,2;0,0,0),(4,0)$),(
%$(1,1,2;1,1,1),(2,0)$),( 
%$(1,2,1;0,1,0),(2,0)$),(
%$(1,2,1;1,0,1),(4,0)$),(
%$(2,1,1;0,1,0),(4,0)$),(
%$(2,1,1;1,0,1),(2,0)$),(
%$(2,2,0;0,0,0),(5,0)$),(
%$(2,0,2;0,0,0),(5,0)$),(
%$(0,2,2;0,0,0),(5,0)$) 
%$\}$

\item \fbox{$(0,0,0)(0,1,0)(1,0,1)\leadsto(37,0)$}\\
$\{$
($(0,0,0;0,0,0),(3,0)$),
($(1,3,0;0,0,0),(2,0)$),
($(0,1,3;0,1,0),(2,0)$),
($(1,1,2;0,0,0),(4,0)$),
($(1,1,2;1,0,1),(2,0)$), 
($(1,2,1;0,1,0),(2,0)$),
($(1,2,1;1,1,1),(2,0)$), 
($(2,1,1;0,1,0),(4,0)$),
($(2,1,1;1,1,1),(2,0)$),
($(2,2,0;0,0,0),(4,0)$),
($(2,0,2;0,0,0),(5,0)$),
($(2,0,2;1,0,1),(1,0)$),
($(0,2,2;0,0,0),(4,0)$)
$\}$

\item \fbox{$(0,0,0)(0,1,0)(1,1,1)\leadsto(42,0)$}\\
$\{$
($(0,0,0;0,0,0),(3,0)$),
($(1,3,0;0,0,0),(2,0)$),
($(1,3,0;1,1,0),(4,0)$),
($(3,0,1;1,0,0),(4,0)$),
($(0,1,3;0,1,0),(4,0)$),
($(1,1,2;0,0,0),(4,0)$),
($(1,1,2;1,1,0),(6,0)$),
($(1,2,1;0,1,0),(2,0)$),
($(1,2,1;1,0,0),(6,0)$),
($(2,1,1;0,1,0),(6,0)$),
($(2,1,1;1,0,0),(2,0)$),
($(2,2,0;0,0,0),(5,0)$),
($(2,2,0;1,1,0),(1,1)$),
($(2,0,2;0,0,0),(6,0)$),
($(0,2,2;0,0,0),(6,0)$)
$\}$

\end{itemize}

%%%%%%%%%%%%%%%%%%%%%%%%%%%%%%%
\subsection{Quotients of $V_6$}
%%%%%%%%%%%%%%%%%%%%%%%%%%%%%%%

For $n=6$, we have a single orbifold. In order to further shorten the notation, we will use $S_3$ symmetry in $(\IZ/6\IZ)^3$: we will write one element per $S_3$ orbit (the size of the orbit is is written between square brackets):

$\{$
$((0,0,0;0,0,0),(3,0))[1],$
$((1,5,0;0,0,0),(6,0))[6],$
$((1,4,1;0,0,0),(6,0))[3],$
$((1,3,2;0,0,0),(24,0))[6],$
$((2,4,0;0,0,0),(24,0))[6],$
$((2,2,2;0,0,0),(5,0))[1],$  
$((3,3,0;0,0,0),(12,0))[3]$
$\}$\\
$(h_{11},h_{12})=(80,0)$.

%%%%%%%%%%%%%%%%%%%%%%%%%%%%%%%
\section{The fundamental Group}
%%%%%%%%%%%%%%%%%%%%%%%%%%%%%%%

We will compute $\pi_1$ of our orbifolds using the fact that they are the quotients of a simply connected space ($\IC^3$).\\ 
Consider $E_n\times E_n\times E_n/G$ as the quotient of $\IC^3$ by $\tilde G$, extension of $G$ by the lattice group $\Lambda_n$.  Let $F=\{g\in \tilde G: \exists x\in \IC^3 \ |\  g.x=x \} $ and $N(F)$ the group generated by $F$.

\begin{theorem} 
The fundamental group of $E_n\times E_n\times E_n/G$ is $\tilde G/N(F)$.
\end{theorem}

A proof can be found in \cite{BroHi} or a more heuristic argument is given in \cite{DHVW2}.

%%%%%%%%%%%%%%%%%%
\subsection{Order 4 and 6}
%%%%%%%%%%%%%%%%%%

\begin{proposition}
Let $G$ be a subgroup of the Vafa-Witten Group with $n=4$ or $6$ which surjects onto the multiplicative part. The closure $N(F)$ of $F$ is the whole of $\tilde G$. 
\end{proposition}

\begin{proof}
We will generalize slightly the notation used up to now, we will write 
$g=(a_1,a_2,a_3;\tau_1 t_i,\tau_2 t_i,\tau_3 t_i)$ where before we would have omitted the $t_i$. This permits to write the elements of the lattice group $\Lambda$ in the form $(0,0,0;\alpha_1 + \beta_1 \omega,\alpha_2 + \beta_2 \omega,\alpha_3 + \beta_3 \omega)$. Let $g$ be an elements of $\tilde G-F$, it is of the form $(a_1,a_2,a_3;v_1,v_2,v_3)$ where, for at least one $k$, $v_k$ is not zero and $a_k=0$. \\
The key is that $\tilde G$ always contains an element of the form $h=(b_1,b_2,b_3;*,*,*)$ such that $a_k+b_k\neq 0$ and $b_k\neq 0$ for all $k\in\{1,2,3\}$. 

{\bf Case 1:} only one of the $a_k=0$. We can assume without loss of generality that $k=1$. Pick $\epsilon,\delta\in\IZ/i$ such that 
$\epsilon,\delta\neq 0$,$\delta\neq -\epsilon$,$\delta\neq a_2$ and $\delta\neq a_2 - \epsilon$. For $n$ large enough (that is $n>3$) we see that there exists $n^2-5n+6>0$ such possible pairs. Indeed the previous conditions take away (in the right order) $n$, $n-1$, $n-1$, $n-2$ and $n-2$ possibilities from the $n^2$ possible pairs of $\IZ/n\times\IZ/n$. (See figure \ref{polycomb})
Now, we just pick $h=(\epsilon,\delta-a_2,a_2-\delta+\epsilon;*,*,*)$.  

\begin{figure}
\centering
\includegraphics[height=3cm]{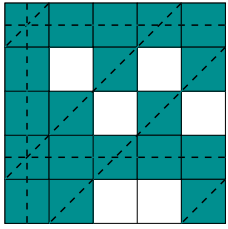}\caption{``Good pairs'' of indices.}
\label{polycomb}
\end{figure}

{\bf Case 2:} all $a_k=0$. We can pick $h$ to be $(1,1,n-2;*,*,*)$. Since $n>3$, $n-2\neq 0$. \\
In both cases $h$ and $h^{-1}.g$ belong to $F$ so from the trivial equality $g=h.(h^{-1}.g)$ we deduce that $g\in N(F)$.  

\end{proof}

\begin{corollary}
Let $G$ be a subgroup of the Vafa-Witten Group which surjects onto the multiplicative part with $n=4$ or $6$, $\pi_1(X(G),x)$ is trivial.
\end{corollary}

%%%%%%%%%%%%%%%%
\subsection{Order 3}
%%%%%%%%%%%%%%%%

\begin{proposition}
All orbifolds obtained when is $n$ is $3$ are simply connected, except the quotient by $III.4$ whose fundamental group has order $3$.
\end{proposition}

\begin{proof}
Let $j=(0,0,0;v_1,v_2,v_3)$ be a translation element in $\tilde G$. We can decompose it as 
$$j=(j g_2 ^{-2} g_1^{-1})(g_1 g_2 ^2).$$
Both $g_1 g_2 ^2=(2,2,2;*,*,*)$ and $j g_2 ^{-2} g_1^{-1}=(1,1,1;*,*,*)$ belong to $F$ and hence $j$ belongs to $F$ as well. The group of translations, $T$, in $\tilde G$ is a normal subgroup (as is the case in any Euclidean group) and so, $G/N(F)\cong (G/T)/(N(F)/T)$. Since $G/T$ is generated by the classes of $g_1$ and $g_2$ and $g_1 g_2^2\in F$, the group $G/N(F)$ is a quotient of $\IZ/3\IZ$. \\
For all groups $III.x$, $x$ different from $4$, $g_1$ is an element of $F$ so the fundamental group is trivial. For $III.4$, it is easy to see that only powers of $g_1 g_2^2$ lie in $F$ and so the fundamental group is cyclic of order $3$.   
\end{proof}

%\begin{appendix}

%<insert any appendices here>

%\section{Appendice}

%\end{appendix}

%\bibliographystyle{my-h-elsevier}

\end{document}